\documentclass[a4paper,11pt]{article}
\usepackage[bookmarks=true,
              bookmarksnumbered=true, breaklinks=true,
              pdfstartview=FitH, hyperfigures=false,
              plainpages=false, naturalnames=true,
              colorlinks=true,
              pdfpagelabels]{hyperref}
\usepackage{geometry,latexsym,amssymb,amsmath,amsthm,color,bm}
\usepackage[latin5]{inputenc}
\usepackage{enumerate}
\usepackage{enumitem}
\usepackage[T1]{fontenc}
\usepackage{authblk}
\usepackage[all]{xy}
\usepackage{palatino}
\usepackage{indentfirst}
\usepackage{titlesec}
\usepackage{graphics}
\usepackage{tabularx}
\usepackage{lipsum}
\usepackage{longtable}
\usepackage{array}
\usepackage{fancyvrb,shortvrb}
\usepackage[ruled,vlined]{algorithm2e}

\newcommand{\GAP}      {{\sf GAP}}
\newcommand{\groupoids}{{\sf groupoids}}
\newcommand{\HAP}      {{\sf HAP}}

\newcommand{\XMod}     {{\sf XMod}}
\newcommand{\XModAlg}  {{\sf XModAlg}}

\newcommand{\catCatt}   {{\bf Cat2}}

\newcommand{\catXSq}    {{\bf XSq}}

\newcommand{\Aut}   {\mathrm{Aut}}

\newcommand{\id}    {\mathrm{id}}

\newcommand{\im}    {\mathrm{im}\,}
\newcommand{\Inn}   {\mathrm{Inn}}

\newcommand{\calC}{\mathcal{C}}

\newcommand{\calS}{\mathcal{S}}

\newcommand{\calX}{\mathcal{X}}

\newcommand{\bt}      {\boxtimes}
\newcommand{\btt}     {\ \tilde{\boxtimes}\ }

\theoremstyle{plain}
\newtheorem{theorem}{Theorem}[section]

\theoremstyle{definition}
\newtheorem{definition}[theorem]{Definition}
\newtheorem{example}[theorem]{Example}

\begin{document}
\title{Computing 3-Dimensional Groups : Crossed Squares and  Cat$^2$-Groups}

\author[a]{Z. Arvasi}
\author[a]{A. Odaba\c{s}}
\author[b]{C.~D. Wensley}
\affil[a]{\small{Department of Mathematics and Computer Science, Osmangazi University, Eskisehir, Turkey}}
\affil[b]{\small{School of Computer Science and Electronic Engineering, Bangor University, North Wales, UK}}

\date{}

\maketitle

\begin{abstract}
The category \catXSq\ of crossed squares is equivalent to 
the category \catCatt\ of cat$^2$-groups. 
Functions for computing with these structures have been developed in 
the package \XMod\ written using the \GAP\ computational discrete algebra 
programming language.
This paper includes details of the algorithms used. 
It contains tables listing the $1,000$ isomorphism classes of cat$^2$-groups on groups of order at most $30$.

\end{abstract}

\noindent{\bf Key Words:} cat$^2$-group, crossed square, \GAP, \XMod\ 
\\ {\bf Classification:} 18D35, 18G50.

\section{Introduction}

This paper is concerned with the latest developments in the general programme 
of "computational higher-dimensional group theory" which forms part of the 
"higher-dimensional group theory" programme described, for example, 
by Brown in \cite{brown-lms}. 

The $2$-dimensional part of these programmes is concerned with group objects 
in the categories of groups or groupoids, and these objects may equivalently 
be considered as crossed modules or cat$^1$-groups. 
A summary of the definitions of these objects, 
with some examples, is contained in \S 2. 

The initial computational part of this programme was described in 
Alp and Wensley \cite{alp-wensley-ijac}. 
The output from this work was the package \XMod\ \cite{xmod} 
for \GAP\ \cite{gap} which, at the time, 
contained functions for constructing crossed modules and cat$^1$-groups of groups, and their morphisms, and conversions from one to another. 
It also contained functions for computing the monoid of derivations of a 
crossed module, and the equivalent monoid of sections of a cat$^1$-group. 
The next development of \XMod\ used the package \groupoids\ 
\cite{groupoids} to compute crossed modules of groupoids. 
Later still, a \GAP\ package \XModAlg\ \cite{xmodalg} 
was written to compute cat$^1$-algebras and crossed modules of algebras, 
as described in \cite{arvasi-odabas}. 

The $3$-dimensional part of the higher-dimensional group theory programme 
is concerned with objects in the category \catXSq\ of crossed squares 
and the equivalent cat$^2$-groups category \catCatt. 
The mathematical basis of these structures is described in \S 3, 
and some computational details are included in \S 4. 
In \S 5 we enumerate the $1,000$ isomorphism 
classes of cat$^2$-group structures on the $92$ groups of order at most $30$. 

There are many other ways of viewing crossed squares and cat$^2$-groups. 
Conduch\'{e} in \cite{conduche} defined the equivalent notion of 
$2$-crossed module.
Brown and Gilbert in \cite{brown-gilbert} introduced braided, 
regular crossed modules as an alternative algebraic model of homotopy $3$-types. 
They also proved that these structures are equivalent to simplicial groups 
with Moore complex of length $2$. 
In \cite{arvasi-ulualan} Arvasi and Ulualan explore the algebraic relationship 
between these structures and also the quadratic modules of Baues \cite{baues}, 
and the homotopy equivalences between them. 

The impetus for the study of higher-dimensional groups 
comes from algebraic topology \cite{brown-indag}. 
Crossed modules are algebraic models of connected (weak homotopy) $2$-types, 
while crossed squares model connected $3$-types. 
The principal topological example of a crossed module arises from 
a pointed pair of spaces $ A \subseteq X$ where the boundary map is 
$\partial : \pi_2(X,A) \to \pi_1(A)$. 
Similarly, given a triad of pointed spaces $A \subseteq X,\ B \subseteq X$ 
we obtain a crossed square as shown in the left-hand diagram below. 
A simple case, when $X$ is a $2$-sphere and $A,B$ are the 
upper and lower hemispheres, results in the square on the right. 
Here $F$ is a free group on one generator $x$, 
the boundaries are the trivial and identity homomorphisms, 
and the crossed pairing is given by 
$h : F \times F \to F,\ (x^i,x^j) \mapsto x^{ij}$  
(see Ellis \cite{ellis}). 
\begin{equation*} 
\xymatrix@R=20pt@C=20pt
{     &  \pi_3(X;A,B) \ar[dd] \ar[rr] \ar[ddrr] 
         &  & \pi_2(B, A \cap B) \ar[dd] 
               &  &  &  F \ar[dd]_{0} \ar[rr]^{0} \ar[ddrr]^{0} 
                        &  &  F \ar[dd]^{{\rm id}} \\
      &  &  &  &  &    \\  
      &  \pi_2(A,A \cap B) \ar[rr] 
         &  & \pi_1(A \cap B) & &  &  F \ar[rr]_{{\rm id}} 
                        &  &  F } 
\end{equation*}

The \XMod\ package follows a purely algebraic approach, 
and does not compute any specifically topological results. 
The interested reader may wish to investigate the \GAP\ package
\HAP\ \cite{hap} which also computes with cat$^1$-groups.

\section{Crossed Modules and Cat$^{1}$-Groups}

The notion of crossed module, generalizing the notion of a G-module, 
was introduced by Whitehead \cite{whitehead-II} in the course of his studies 
on the algebraic structure of the second relative homotopy group.

A \emph{crossed module} consists of a group homomorphism 
$\partial : S \rightarrow R$, endowed with a left action $R$ on $S$ 
(written by $(r,s) \rightarrow {}^{r}s$ for $r \in R$ and $s \in S$) 
satisfying the following conditions:

\begin{center}
	\begin{tabular}{rclll}
	$\partial (^{r}s)$ 
		& $=$ 
			& $r(\partial s)r^{-1}$ 
				&   & $\forall~ s \in S,~ r \in R; $ \\
	$^{(\partial s_{2})}s_{1}$ 
		& $=$ 
			& $s_{2}s_{1}s_{2}^{-1}$ 
				&   & $\forall~ s_{1},s_{2} \in S$. 
	\end{tabular}
\end{center}

The first condition is called the \emph{pre-crossed module property} 
and the second one the \emph{Peiffer identity}. 
We will denote such a crossed module by $\calX = (\partial : S \rightarrow R)$.

A \emph{morphism of crossed modules} 
$(\sigma ,\rho ) : \calX_{1} \rightarrow \calX_{2}$, 
where $\calX_{1} = (\partial_{1} : S_{1} \rightarrow R_{1})$ 
and   $\calX_{2} = (\partial_{2} : S_{2} \rightarrow R_{2})$, 
consists of two group homomorphisms $\sigma : S_{1} \rightarrow S_{2}$
and $\rho : R_{1} \rightarrow R_{2}$ such that 
\[ 
\partial_{2}\circ\sigma ~=~ \rho\circ\partial_{1}, 
\quad \mbox{and} \quad 
\sigma(^{r}s) ~=~ ^{(\rho r)}\sigma s 
\qquad
\forall s \in S,~ r \in R.
\] 

Standard constructions for crossed modules include the following. 
\begin{enumerate}
\item 
A \emph{conjugation crossed module} \index{conjugation crossed module} 
is an inclusion of a normal subgroup $N \unlhd R$, 
where $R$ acts on $N$ by conjugation.
\item 
An \emph{automorphism crossed module} \index{automorphism crossed module} 
has as range a subgroup $R$ of the automorphism group $\Aut(S)$ of $S$ 
which contains the inner automorphism group $\Inn(S)$ of $S$. 
The boundary maps $s \in S$ to the inner automorphism of $S$ by $s$.
\item 
A \emph{zero boundary crossed module} \index{$R$-module} 
has a $R$-module as source and $\partial = 0$.
\item 
Any homomorphism $\partial : S \to R$, with $S$ abelian 
and $\im\partial$ in the centre of $R$, 
provides a crossed module with $R$ acting trivially on $S$.
\item 
A \emph{central extension crossed module} 
\index{central extension crossed module} 
has as boundary a surjection $\partial : S \to R$ with central kernel, 
where $r \in R$ acts on $S$ by conjugation with $\partial^{-1}r$.
\item 
The \emph{direct product} of \index{direct product!of crossed modules} 
$\calX_1 = (\partial_1 : S_1 \to R_1)$ and $\calX_2 = (\partial_2 : S_2 \to R_2)$ 
is $\calX_1 \times \calX_2 
= (\partial_1 \times \partial_2 : S_1 \times S_2 \to R_1 \times R_2)$ 
with direct product action 
${}^{(r_1,r_2)}(s_1,s_2) = \left({}^{r_1}s_1,{}^{r_2}s_2\right)$. 
\end{enumerate}

Loday reformulated the notion of crossed module as a cat$^{1}$-group. 
Recall from \cite{Loday} that a \emph{cat$^{1}$-group} is a triple $(G;t,h)$ consisting of a group $G$ with two endomorphisms: 
the \emph{tail map} $t$ and the \emph{head map} $h$, 
having a common image $R$ and satisfying the following axioms. 
\begin{equation} \label{cat1-axioms} 
t \circ h = h, \quad  
h \circ t = t, 
\quad \mbox{and}\quad  [\ker t,\ker h] = 1. 
\end{equation} 
When only the first two of these axioms are satisfied, the structure is a 
\emph{pre-cat$^1$-group}. 
It follows immediately that $t \circ t = t$ and $h \circ h = h$. 
We picture $(G;t,h)$ as 
\[
\xymatrix@R=40pt@C=40pt
{ G \ar@<+0.4ex>[r]^{t,h} \ar@<-0.4ex>[r]  &  R }
\] 

A \emph{morphism of cat}$^{1}$\emph{-groups} 
$(G_{1};t_1,h_1) \rightarrow (G_{2};t_2,h_2)$ 
is a group homomorphism $f : G_{1} \rightarrow G_{2}$ such that 
\[ 
f \circ t_1 ~=~ t_2 \circ f  
\quad\mbox{and}\quad 
f \circ h_1 ~=~ h_2 \circ f.
\] 
Crossed modules and cat$^{1}$-groups are equivalent two-dimensional 
generalisations of a group. 
It was shown in \cite[Lemma 2.2]{Loday} that, 
on setting $S = \ker t,~ R = \im t$ and $\partial = h|_{S}$, 
the conjugation action makes $(\partial : S \rightarrow R)$ 
into a crossed module. 
Conversely, if $(\partial : S \rightarrow R)$ is a crossed module, 
then setting $G = S \rtimes R$ and defining $t,h$ by $t(s,r) = (1,r)$ 
and $h(s,r) = (1,(\partial s)r)$ for $s \in S$, $r \in R$, 
produces a cat$^{1}$-group $(G;t,h)$.

\section{Crossed Squares and Cat$^{2}$-Groups}

The notion of a crossed square is due to Guin-Walery and Loday \cite{walery-loday}. 
A \emph{crossed square of groups} $\calS$ is a commutative square of groups 
\begin{equation} \label{xsq-diag}
\xymatrix@R=20pt@C=20pt
{     &  L \ar[dd]_{\lambda} \ar[rr]^{\kappa} \ar[ddrr]^{\pi} 
         &  & M \ar[dd]^{\mu} 
               &  &  &  L \ar[dd]_{\kappa} \ar[rr]^{\lambda} \ar[ddrr]^{\pi} 
                        &  &  N \ar[dd]^{\nu} \\
\calS \quad = 
      &  &  &  &  &  \tilde{\calS} \quad =  \\  
      &  N \ar[rr]_{\nu} 
         &  & P & &  &  M \ar[rr]_{\mu} 
                        &  &  P } 
\end{equation}
\noindent together with left actions of $P$ on $L,M,N$ 
and a \emph{crossed pairing} map ${\ \bt\ } : M \times N \rightarrow L$. 
Then $M$ acts on $N$ and $L$ via $P$ and $N$ acts on $M$ and $L$ via $P$. 
The diagram illustrates an \emph{oriented crossed square}, 
since a choice of where to place $M$ and $N$ has been made. 
The \emph{transpose} $\tilde{\calS}$ of $\calS$ is obtained by making the alternative choice. 
Since crossed pairing identities are similar to those for commutators, 
the crossed pairing for $\tilde{\calS}$ is $\btt$ 
where $(n \btt m) = (m \bt n)^{-1}$. 
Transposition gives an equivalence relation on the set of 
oriented crossed squares, and a crossed square is an equivalence class. 

The structure of an oriented crossed square must satisfy the following axioms 
for all $l \in L,~ m,m^{\prime} \in M,~ n,n^{\prime} \in N$ and $p \in P$. 
\begin{enumerate} 
\item 
With the given actions, the homomorphisms $\kappa, \lambda, \mu, \nu$ 
and $\pi = \mu\circ\kappa = \nu\circ\lambda$ are crossed modules, 
and both $\kappa, \lambda$ are $P$-equivariant, 
\item
$(mm^{\prime} {\ \bt\ } n) ~=~ (^{m}m^{\prime} {\ \bt\ } {^{m}n})\,(m {\ \bt\ } n)$  
\quad and \quad 
$(m {\ \bt\ } {nn^{\prime}}) ~=~ (m {\ \bt\ } n)\,(^{n}m {\ \bt\ } ^{n}n^{\prime})$, 
\item 
$\kappa(m {\ \bt\ } n) ~=~ m({}^{n} m^{-1})$ 
\quad and \quad 
$\lambda(m {\ \bt\ } n) ~=~ ({}^{m}n)n^{-1}$, 
\item 
$(\kappa l {\ \bt\ } n) ~=~ l({}^{n} l^{-1})$  
\quad and \quad 
$(m {\ \bt\ } \lambda l) ~=~ ({}^{m}l)l^{-1}$, 
\item 
$^{p}(m {\ \bt\ } n) ~=~ (^{p}m {\ \bt\ } ^{p}n)$. 
\end{enumerate} 

\noindent 
Note that axiom 1. implies that $(\id,\mu), (\id,\nu), (\kappa,\id)$ 
and $(\lambda,\id)$ are morphisms of crossed modules. 

\medskip 
\noindent 
Standard constructions for crossed squares include the following.
\begin{enumerate}
\item 
If $M,N$ are normal subgroups of the group $P$ then the diagram of inclusions
\[
\xymatrix@R=40pt@C=40pt
{ M \cap N \ar[r]^(0.6){} \ar[d]_{}  
	& M \ar[d]^{} \\
	N \ar[r]_{}  
	& P }
\] 
\noindent together with the actions of $P$ on $M,N$ and $M\cap N$ 
given by conjugation, and the commutator map 
\[
[~,~] ~:~ M\times N \rightarrow M\cap N,\quad 
(m,n)\mapsto [m,n] \,=\, mnm^{-1}n^{-1}, 
\] 
is a crossed square. 
We call this an \emph{inclusion crossed square}.
\item 
The diagram
\[
\xymatrix@R=40pt@C=40pt
{ M \ar[r]^{\alpha} \ar[d]_{\alpha} 
	& \Inn\,M \ar[d]^{\iota} \\
	\Inn\,M \ar[r]_{\iota} 
	& \Aut\,M }
\] 
\noindent is a crossed square, 
where $\alpha $ maps $m\in M$ to the inner automorphism%
\[
\beta_{m} : M \rightarrow M,\quad 
m^{\prime}\mapsto mm^{\prime}m^{-1}, 
\]
and where $\iota $ is the inclusion of $\Inn\,M$ in $\Aut\,M$; 
the actions are standard; and the crossed pairing is
\[
\bt ~:~ \Inn\,M \times \Inn\,M \rightarrow M,\quad 
(\beta_{m},\beta_{m^{\prime}}) \;\mapsto\; [m,m^{\prime}]~.
\]
\item 
If $P$ is a group and $M,N$ are ordinary $P$-modules, 
and if $A$ is an arbitrary abelian group on which $P$ is assumed to act trivially, 
then there is a crossed square
\[
\xymatrix@R=40pt@C=40pt
{ A \ar[r]^{0} \ar[d]_{0}  
	& M \ar[d]^{0} \\
	N \ar[r]_{0} 
	& P }
\]
\item 
Given two crossed modules, $(\mu : M \rightarrow P)$ and $(\nu : N \rightarrow P)$, 
there is a universal crossed square 
\[
\xymatrix@R=40pt@C=40pt
{ M \otimes N \ar[d]_{\lambda} \ar[r]^{\kappa} 
	& M \ar[d]^{\mu} \\ 
	N \ar[r]_{\nu} 
	& P } 
\] 
where $M \otimes N$ is constructed using the nonabelian tensor product of groups 
\cite{brown-loday}. 
\item
The \emph{direct product} of crossed squares $\calS_1,\calS_2$ is 
\[
\xymatrix@R=40pt@C=50pt
{ L_1 \times L_2 \ar[r]^{\kappa_1 \times \kappa_2} 
                 \ar[d]_{\lambda_1 \times \lambda_2}  
	& M_1 \times M_2 \ar[d]^{\mu_1 \times \mu_2} \\
	N_1 \times N_2 \ar[r]_{\nu_1 \times \nu_2} 
	& P_1 \times P_2 }
\] 
with actions 
\[
{}^{(p_1,p_2)}(l_1,l_2) = \left({}^{p_1}l_1,{}^{p_2}l_2\right), \quad 
{}^{(p_1,p_2)}(m_1,m_2) = \left({}^{p_1}m_1,{}^{p_2}m_2\right), \quad 
{}^{(p_1,p_2)}(n_1,n_2) = \left({}^{p_1}n_1,{}^{p_2}n_2\right),  
\]
and crossed pairing 
\[
\bt\left((m_1,m_2),(n_1,n_2)\right) ~=~ \left(\bt_1(m_1,n_1),\bt_2(m_2,n_2)\right). 
\]
\end{enumerate}

The crossed square $\calS$ in (\ref{xsq-diag}) can be thought of 
as a horizontal or vertical crossed module of crossed modules:
\[
\xymatrix@R=20pt@C=20pt
{ L \ar[dd]_{\lambda}  
	&  &  M \ar[dd]^{\mu} 
	      &  &  L \ar[rr]^{\kappa}
	            &  {} \ar[dd]^{(\lambda,\mu)} 
	               &  M \\ 
\quad \ar[rr]^{(\kappa,\nu)} 
    &  &  &  &  &  & \quad \\
  N &  &  P 
	      &  &  N \ar[rr]_{\nu} 
	            &  {} 
	               &  P 
} 
\]
\noindent 
where $(\kappa,\nu)$ is the boundary of the crossed module with 
domain $(\lambda : L \rightarrow N)$ and codomain $(\mu : M \rightarrow P)$. 
(See also section 9.2 of \cite{wensley-notes}.)

There is an evident notion of morphism of crossed squares  
which preserves all the structure, 
so that we obtain a category \catXSq, the category of crossed squares.

\medskip
Although, when first introduced by Loday and Walery \cite{walery-loday}, 
the notion of crossed square of groups was not linked to that of cat$^{2}$-groups, 
it was in this form that Loday gave their generalisation 
to an $n$-fold structure, cat$^{n}$-groups (see \cite{Loday}). 
When $n=1$ this is the notion of cat$^1$-group given earlier.

When $n=2$ we obtain a cat$^{2}$-group. 
Again we have a group $G$, but this time with two \emph{independent} 
cat$^{1}$-group structures on it. 
So a \emph{cat$^{2}$-group} is a $5$-tuple, $(G;t_1,h_1;t_2,h_2)$, 
where $(G;t_{i},h_{i}),~ i=1,2$ are cat$^{1}$-groups, and
\[
t_{1} \circ t_{2} ~=~ t_{2} \circ t_{1}, \quad 
h_{1} \circ h_{2} ~=~ h_{2} \circ h_{1}, \quad 
t_{1} \circ h_{2} ~=~ h_{2} \circ t_{1}, \quad
t_{2} \circ h_{1} ~=~ h_{1} \circ t_{2}. 
\]
To emphasise the relationship with crossed squares 
we may illustrate an \emph{oriented} cat$^2$-group by the diagram 
\begin{equation} \label{cat2-diag}
\xymatrix@R=50pt@C=50pt 
{ G \ar@<+0.4ex>[r]^{t_1,h_1} \ar@<-0.4ex>[r] 
    \ar@<-0.4ex>[d]_{t_2,h_2} \ar@<+0.4ex>[d] 
    \ar@<-0.4ex>[rd]^{\ t_1 \circ t_2} \ar@<+0.4ex>[rd]_{h_1 \circ h_2} 
	& R_1 \ar@<+0.4ex>[d]^{t_2,h_2} \ar@<-0.4ex>[d] \\
  R_2 \ar@<-0.4ex>[r]_{t_1,h_1} \ar@<+0.4ex>[r] 
	& R_{12} 
}
\end{equation} 
where $R_{12}$ is the image of $t_1 \circ t_2 = t_2 \circ t_1$. 

\medskip
A morphism of cat$^{2}$-groups is a triple $(\gamma ,\rho_1 ,\rho_2)$,
as shown in the diagram 
\[
\xymatrix@R=40pt@C=40pt 
{ R_1 \ar[d]_{\rho_1} 
	& G \ar[d]_{\gamma} \ar@<-0.4ex>[l]_{t_1,h_1} \ar@<+0.4ex>[l]  
	                    \ar@<+0.4ex>[r]^{t_2,h_2} \ar@<-0.4ex>[r] 
		& R_2 \ar[d]^{\rho_2} \\
  R'_1 
	& G' \ar@<+0.4ex>[l]^{t'_1,h'_1} \ar@<-0.4ex>[l] 
	     \ar@<-0.4ex>[r]_{t'_2,h'_2} \ar@<+0.4ex>[r] 
		& R'_2 
}
\]
\noindent where 
$\gamma : G \to G^{\prime},~ \rho_1 = \gamma|_{R_1}$ 
and $\rho_2 = \gamma|_{R_2}$ are homomorphisms satisfying: 
\[ 
\rho_1 \circ t_1 = t_1^{\prime} \circ \gamma, \qquad 
\rho_1 \circ h_1 = h_1^{\prime} \circ \gamma, \qquad 
\rho_2 \circ t_2 = t_2^{\prime} \circ \gamma, \qquad 
\rho_2 \circ h_2 = h_2^{\prime} \circ \gamma. 
\] 
We thus obtain a category \catCatt, the category of cat$^{2}$-groups. 

Notice that, unlike the situation with crossed squares 
where the diagonal is a crossed module, 
it is \emph{not} required that the diagonal in (\ref{cat2-diag}) 
is a cat$^1$-group -- it may just be a pre-cat$^1$-group. 
The simplest example of this situation  
is described in Example \ref{ex-d8} below. 

Loday, in \cite{Loday} proved that there is an equivalence of categories 
between the category \catCatt\ and the category \catXSq.
We now consider the sketch proof of this result 
(see also \cite{mutlu-porter-2003}). 

The cat$^{2}$-group $(G;t_1,h_1;t_2,h_2)$ determines a diagram of homomorphisms 
\begin{equation} \label{ker-im-diag}
\xymatrix@R=50pt@C=50pt
{ \ker t_1 \cap \ker t_2 \ar[d]_{(\id,\partial_2)} \ar[r]^{(\partial_1,\id)} 
  	  & \im t_1 \cap \ker t_2 \ar[d]^{(\id,\partial_2)} \\ 
  \ker t_1 \cap \im t_2 \ar[r]_{(\partial_1,\id)}  
	  & \im t_1 \cap \im t_2 } 
\end{equation} 
\noindent where each morphism is a crossed module for the natural action, 
conjugation in $G$. 
The required crossed pairing is given by the commutator in $G$ since, 
if $x \in \im t_1 \cap \ker t_2$ and $y \in \ker t_1 \cap \im t_2$ 
then $[x,y] \in \ker t_1 \cap \ker t_2$. 
It is routine to check the crossed square axioms.
	
Conversely, if
\[
\xymatrix@R=40pt@C=40pt
{ L \ar[d]_{\lambda} \ar[r]^{\kappa}  
	  & M \ar[d]^{\mu} \\ 
  N \ar[r]_{\nu} 
	  & P }  
\]
\noindent is a crossed square, 
then we consider it as a morphism of crossed modules 
$(\kappa,\nu) : (\lambda : L \to N) \rightarrow (\mu : M  \to P)$.
Using the equivalence between crossed modules and cat$^{1}$-groups this
gives a morphism
\[
\partial : (L \rtimes N,t,h) \longrightarrow (M \rtimes P, t^{\prime}, h^{\prime})
\]
of cat$^{1}$-groups. 
There is an action of $(m,p) \in M \rtimes P$ on $(l,n) \in L \rtimes N$ 
given by
\[
^{(m,p)}(l,n) ~=~ (^{m}(^{p}l) (m \bt\ ^{p}n),\ ^{p}n)\,.
\] 
Using this action, we form its associated cat$^{2}$-group with source  
$(L \rtimes N) \rtimes (M \rtimes P)$ 
and induced endomorphisms $t_1,h_1,t_2,h_2$. 

\begin{example} \label{ex-d8} 
Let $D_8 = \langle a,b ~|~ a^2, b^2, (ab)^4 \rangle$ 
be the dihedral group of order $8$, 
and let $c=[a,b]=(ab)^2$ so that $a^b=ac$ and $b^a=bc$.  
(The standard permutation representation is given by 
$a=(1,2)(3,4), b=(1,3), ab=(1,2,3,4), c=(1,3)(2,4)$.) 

Define $t_a,t_b : D_8 \to D_8$ by $t_a : [a,b] \mapsto [a,1]$ 
and $t_b : [a,b] \mapsto [1,b]$. 
Construct cat$^1$-groups $C_a = (D_8,t_a,t_a)$ and $C_b = (D_8,t_b,t_b)$. 
Diagrams (\ref{cat2-diag}) and (\ref{ker-im-diag}) become 
\[
\xymatrix@R=40pt@C=50pt
{ D_8 \ar@<+0.4ex>[r]^{t_a} \ar@<-0.4ex>[r] 
    \ar@<-0.4ex>[d]_{t_b} \ar@<+0.4ex>[d] 
    \ar@<-0.4ex>[rd]^{\ t} \ar@<+0.4ex>[rd] 
	& A \ar@<+0.4ex>[d]^{t_b} \ar@<-0.4ex>[d] 
	   &  & C \ar[r]^{c\ \mapsto\ 1} \ar[d]_{c\ \mapsto\ 1}  
	         & A \ar[d]^{a\ \mapsto\ 1} \\
  B \ar@<-0.4ex>[r]_{t_a} \ar@<+0.4ex>[r] 
	& I 
	   &  & B \ar[r]_{b\ \mapsto\ 1}  
	         & I }
\] 
where $A = \langle a \rangle$, $B = \langle b \rangle$, $C = \langle c \rangle$ 
and $I$ is the trivial group. 
The crossed pairing is given by $\bt(a,b)=c$. 
The diagonal map $t = t_a \circ t_b$ has kernel $D_8$, 
and $[\ker t,\ker t] = C$, so the diagonal is \emph{not} a cat$^1$-group. 
\end{example} 

\begin{definition}
A \emph{cat$^{n}$-group} consists of a group $G$ 
with $n$ independent cat$^{1}$-group structures $(G;t_{i},h_{i})$, 
$1 \leq i \leq n$, such that for all $i \ne j$ 
\[
t_{i}t_{j} = t_{j}t_{i}, \quad 
h_{i}h_{j} = h_{j}h_{i} \quad \mbox{and} \quad 
t_{i}h_{j} = h_{j}t_{i}. 
\]
\end{definition}

A generalisation of crossed square to higher dimensions was given by Ellis
and Stenier (cf. \cite{ellis-stenier}). 
It is called a \textquotedblleft crossed $n$-cube\textquotedblright. 
We only use this construction for $n=2$.

\section{Computer Implementation}

\GAP\ \cite{gap} is an open-source system for discrete computational
algebra. The system consists of a library of implementations of mathematical
structures: groups, vector spaces, modules, algebras, graphs, codes,
designs, etc.; plus databases of groups of small order, character tables, etc. 
The system has world-wide usage in the area of education and scientific research. 
\GAP\ is free software and user contributions to the system are supported. 
These contributions are organized in a form of \GAP\ packages 
and are distributed together with the system.  Contributors can
submit additional packages for inclusion after a reviewing process.

The Small Groups library by Besche, Eick and O'Brien in \cite{besche-eick-obrien} 
provides access to descriptions of the groups of small order. 
The groups are listed up to isomorphism. 
The library contains all groups of order at most 2000 except 1024.

\subsection{2-Dimensional Groups}

The \XMod\ package for \GAP\ contains functions for computing with 
crossed modules, cat$^{1}$-groups and their morphisms, 
and was first described in \cite{xmod}. 
A more general notion of cat$^1$-group is implemented in \XMod, 
where the tail and head maps are no longer required to be endomorphisms on $G$. 
Instead it is required that $t$ and $h$ have a common image $R$, 
and an \emph{embedding} $e : R \to G$ is added.  
The axioms in (\ref{cat1-axioms}) then become:  
\[ 
t \circ e \circ h = h, \quad  
h \circ e \circ t = t, 
\quad \mbox{and}\quad  [\ker t,\ker h] = 1, 
\] 
and again it follows that $t \circ e \circ t = t$ and $h \circ e \circ h = h$. 
We denote such a cat$^1$-group by $(e;t,h : G \to R)$.

This package may be used to select a cat$^{1}$-group from a data file. 
All cat$^{1}$-structures on groups of size up to 70 
(ordered according to the \GAP\ numbering of small groups) 
are stored in a list in the file \texttt{cat1data.g}.
The function \textbf{Cat1Select} may be used in three ways. 
\textbf{Cat1Select( size )} returns the names of the groups with this size, 
while \textbf{Cat1Select( size, gpnum )} prints a list of cat$^1$1-structures 
for this chosen group. 
\textbf{Cat1Select( size, gpnum, num )} returns the chosen cat$^1$1-group.
\textbf{XModOfCat1Group} produces the associated crossed module. 

The following \GAP\ session illustrates the use of these functions.

\begin{Verbatim}[frame=single, fontsize=\small, commandchars=\\\{\}]
\textcolor{blue}{gap> Cat1Select( 12 );}
Usage:  Cat1Select( size, gpnum, num );
[ "C3 : C4", "C12", "A4", "D12", "C6 x C2" ]
\textcolor{blue}{gap> Cat1Select( 12, 3 );}
Usage:  Cat1Select( size, gpnum, num );
There are 2 cat1-structures for the group A4.
Using small generating set [ f1, f2 ] for source of homs.
[ [range gens], [tail genimages], [head genimages] ] :-
(1)  [ [ f1 ], [ f1, <identity> of ... ], [ f1, <identity> of ... ] ]
(2)  [ [ f1, f2 ],  tail = head = identity mapping ]
2
\textcolor{blue}{gap> C1 := Cat1Select( 12, 3, 2 );}
[A4=>A4]
\textcolor{blue}{gap> X1 := XModOfCat1Group( C1 );}
[triv->A4]
\end{Verbatim}

\subsection{3-dimensional Groups}

We have developed new functions for \XMod\ which construct 
(pre-)cat$^{2}$-groups, (pre-)cat$^{3}$-groups, and their morphisms. 
Functions for (pre-)cat$^{2} $-groups include \textbf{PreCat2Group},
\textbf{Cat2Group}, \textbf{IsPreCat2Group}, \textbf{IsCat2Group} 
and \textbf{PreCat2GroupByPreCat1Groups}. 
Attributes of a (pre)cat$^{2}$-group constructed in this way include 
\textbf{GeneratingCat1Groups}, \textbf{Size}, \textbf{Name} and 
\textbf{Edge2DimensionalGroup} where '\textbf{Edge}' is one of 
\{\textbf{Up, Left, Right, Down, Diagonal}\}.  

As with cat$^1$-groups, we use a more general notion for cat$^2$-groups. 
An \emph{oriented cat$^2$-group} has the form 
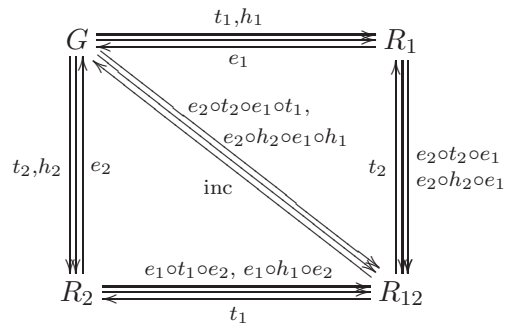
\[
\xymatrix@R=80pt@C=100pt 
{ G \ar@<+0.5ex>[r]^{t_1,h_1} \ar@<+0.1ex>[r] 
    \ar@<-0.5ex>[d]_{t_2,h_2} \ar@<-0.1ex>[d] 
    \ar@<+0.5ex>[rd]^(0.35){e_2 \circ t_2 \circ e_1 \circ t_1,} 
    \ar@<+0.1ex>[rd]^(0.45){\ e_2 \circ h_2 \circ e_1 \circ h_1} 
	& R_1 \ar@<+0.5ex>[d]^(0.45){e_2 \circ t_2 \circ e_1} 
	      \ar@<+0.5ex>[d]^(0.55){e_2 \circ h_2 \circ e_1} 
	      \ar@<+0.1ex>[d] 
	      \ar@<+0.4ex>[l]^{e_1} \\
  R_2 \ar@<+0.5ex>[r]^{e_1 \circ t_1 \circ e_2,~ e_1 \circ h_1 \circ e_2} 
      \ar@<+0.1ex>[r] 
      \ar@<-0.4ex>[u]_{e_2} 
	& R_{12} \ar@<+0.4ex>[l]^{t_1} \ar@<+0.4ex>[u]^{t_2} 
	         \ar@<+0.4ex>[ul]^{{\rm inc}} 
}
\]
where $R_1, R_2$ need not be subgroups of $G$, 
but $R_{12}$ is taken to be the common image of 
$e_2 \circ t_2 \circ e_1 \circ t_1$ and 
$e_1 \circ t_1 \circ e_2 \circ t_2$, a subgroup of $G$.

Generalizing these functions, we have introduced \textbf{Cat3Group} and \textbf{HigherDimension} which construct cat$^{3}$-groups. 
Functions for cat$^{n}$-groups of higher dimension will be added as time permits. 
An orientation of a cat$^{3}$-group on $G$ displays a cube whose six faces 
(ordered as front; up, left, right, down, back) are cat$^{2}$-groups. 
The group $G$ is positioned where the front, up and left faces meet. 
The following \GAP\ session illustrates the use of these functions. 
Notice that the cat$^2$-group \verb|C2ab| is the second example 
with a diagonal which is only a pre-cat$^1$-group. 

\begin{Verbatim}[frame=single, fontsize=\small, commandchars=\\\{\}]
\textcolor{blue}{gap> a := (1,2,3,4)(5,6,7,8);;}
\textcolor{blue}{gap> b := (1,5)(2,6)(3,7)(4,8);;}
\textcolor{blue}{gap> c := (2,6)(4,8);;}
\textcolor{blue}{gap> G := Group( a, b, c );;}
\textcolor{blue}{gap> SetName( G, "c4c2:c2" );}
\textcolor{blue}{gap> t1a := GroupHomomorphismByImages( G, G, [a,b,c], [(),(),c] );; }
\textcolor{blue}{gap> C1a := PreCat1GroupByEndomorphisms( t1a, t1a );;}
\textcolor{blue}{gap> t1b := GroupHomomorphismByImages( G, G, [a,b,c], [a,(),()] );;}
\textcolor{blue}{gap> C1b := PreCat1GroupByEndomorphisms( t1b, t1b );;}
\textcolor{blue}{gap> C2ab := Cat2Group( C1a, C1b );}
(pre-)cat2-group with generating (pre-)cat1-groups:
1 : [c4c2:c2 => Group( [ (), (), (2,6)(4,8) ] )]
2 : [c4c2:c2 => Group( [ (1,2,3,4)(5,6,7,8), (), () ] )]
\textcolor{blue}{gap> IsCat2Group( C2ab );}
true
\textcolor{blue}{gap> Size( C2ab );}
[ 16, 2, 4, 1 ]
\textcolor{blue}{gap> IsCat1Group( Diagonal2DimensionalGroup( C2ab ) );} 
false
\textcolor{blue}{gap> t1c := GroupHomomorphismByImages( G, G, [a,b,c], [a,b,c] );;}
\textcolor{blue}{gap> C1c := PreCat1GroupByEndomorphisms( t1c, t1c );;}
\textcolor{blue}{gap> C3abc := Cat3Group( C1a, C1b, C1c );}
(pre-)cat3-group with generating (pre-)cat1-groups:
1 : [c4c2:c2 => Group( [ (), (), (2,6)(4,8) ] )]
2 : [c4c2:c2 => Group( [ (1,2,3,4)(5,6,7,8), (), () ] )]
3 : [c4c2:c2 => Group( [ (1,2,3,4)(5,6,7,8), (1,5)(2,6)(3,7)(4,8),
(2,6)(4,8) ] )]
\textcolor{blue}{gap> IsPreCat3Group( C3abc );}
true
\textcolor{blue}{gap> HigherDimension( C3abc );}
4
\textcolor{blue}{gap> Front3DimensionalGroup( C3abc );} 
(pre-)cat2-group with generating (pre-)cat1-groups:
1 : [c4c2:c2 => Group( [ (), (), (2,6)(4,8) ] )]
2 : [c4c2:c2 => Group( [ (1,2,3,4)(5,6,7,8), (), () ] )]
\end{Verbatim}

\subsection{Morphisms of 3-Dimensional Groups}

The function \textbf{MakeHigherDimensionalGroupMorphism} defines morphisms of 
higher dimensional groups, such as cat$^{2}$-groups and crossed squares. 
Functions for cat$^{2}$-group morphisms include 
\textbf{Cat2GroupMorphismByCat1GroupMorphisms}, \textbf{Cat2GroupMorphism} and 
\textbf{IsCat2GroupMorphism}. 
The function \textbf{AllCat2GroupMorphisms} is used to find 
all morphisms between two cat$^{2}$-groups.

In the following \GAP\ session, we obtain a cat$^{2}$-group morphism
using these functions.

\begin{Verbatim}[frame=single, fontsize=\small, commandchars=\\\{\}]
\textcolor{blue}{gap> C2_82 := Cat2Group( Cat1Group(8,2,1), Cat1Group(8,2,2) );}
(pre-)cat2-group with generating (pre-)cat1-groups:
1 : [C4 x C2 => Group( [ <identity> of ..., <identity> of ...,
<identity> of ... ] )]
2 : [C4 x C2 => Group( [ <identity> of ..., f2 ] )]
\textcolor{blue}{gap> C2_83 := Cat2Group( Cat1Group(8,3,2), Cat1Group(8,3,3) );}
(pre-)cat2-group with generating (pre-)cat1-groups:
1 : [D8 => Group( [ f1, f1 ] )]
2 : [D8=>D8]
\textcolor{blue}{gap> up1 := GeneratingCat1Groups( C2_82 )[1];;}
\textcolor{blue}{gap> lt1 := GeneratingCat1Groups( C2_82 )[2];;}
\textcolor{blue}{gap> up2 := GeneratingCat1Groups( C2_83 )[1];;}
\textcolor{blue}{gap> lt2 := GeneratingCat1Groups( C2_83 )[2];;}
\textcolor{blue}{gap> G1 := Source( up1 );; R1 := Range( up1 );; Q1 := Range( lt1 );;}
\textcolor{blue}{gap> G2 := Source( up2 );; R2 := Range( up2 );; Q2 := Range( lt2 );;}
\textcolor{blue}{gap> homG := AllHomomorphisms( G1, G2 );;}
\textcolor{blue}{gap> homR := AllHomomorphisms( R1, R2 );;}
\textcolor{blue}{gap> homQ := AllHomomorphisms( Q1, Q2 );;}
\textcolor{blue}{gap> upmor := Cat1GroupMorphism( up1, up2, homG[1], homR[1] );;}
\textcolor{blue}{gap> IsCat1GroupMorphism( upmor );}
true
\textcolor{blue}{gap> ltmor := Cat1GroupMorphism( lt1, lt2, homG[1], homQ[1] );;}
\textcolor{blue}{gap> mor2 := PreCat2GroupMorphism( C2_82, C2_83, upmor, ltmor );}
<mapping: (pre-)cat2-group with generating (pre-)cat1-groups:
1 : [C4 x C2 => Group( [ <identity> of ..., <identity> of ..., 
  <identity> of ... ] )]
2 : [C4 x C2 => Group( [ <identity> of ..., f2 ] )] -> (pre-)cat
2-group with generating (pre-)cat1-groups:
1 : [D8 => Group( [ f1, f1 ] )]
2 : [D8=>D8] >
\textcolor{blue}{gap> IsCat2GroupMorphism( mor2 );}
true
\textcolor{blue}{gap> mor8283 := AllCat2GroupMorphisms( C2_82, C2_83 );;}
\textcolor{blue}{gap> Length( mor8283 );}
2
\end{Verbatim}

\subsection{Natural Equivalence}

The equivalence between categories \catXSq\ and \catCatt\ 
is implemented by the functions 
\textbf{CrossedSquareOfCat2Group} and \textbf{Cat2GroupOfCrossedSquare} 
which construct crossed squares and cat$^{2}$-groups 
from the given cat$^{2}$-groups and crossed squares, respectively.

The following \GAP\ session illustrates the use of these functions.
The dihedral group $D_{20}$ has two normal subgroups $D_{10}$ 
whose intersection is the cyclic $C_5$.  
We construct the crossed square of normal subgroups, 
and then use the conversion functions to obtain the associated cat$^{2}$-group. 
We then obtain the crossed square \texttt{Xab} 
associated to the cat$^2$-group \texttt{C2ab} obtained earlier. 

\begin{Verbatim}[frame=single, fontsize=\small, commandchars=\\\{\}]
\textcolor{blue}{gap> d20 := DihedralGroup( IsPermGroup, 20 );;}
\textcolor{blue}{gap> gend20 := GeneratorsOfGroup( d20 ); }
[ (1,2,3,4,5,6,7,8,9,10), (2,10)(3,9)(4,8)(5,7) ]
\textcolor{blue}{gap> p1 := gend20[1];;  p2 := gend20[2];;  p12 := p1*p2; }
(1,10)(2,9)(3,8)(4,7)(5,6)
\textcolor{blue}{gap> d10a := Subgroup( d20, [ p1^2, p2 ] );; }
\textcolor{blue}{gap> d10b := Subgroup( d20, [ p1^2, p12 ] );; }
\textcolor{blue}{gap> c5d := Subgroup( d20, [ p1^2 ] );; }
\textcolor{blue}{gap> SetName( d20, "d20" );  SetName( d10a, "d10a" ); }
\textcolor{blue}{gap> SetName( d10b, "d10b" );  SetName( c5d, "c5d" );  }
\textcolor{blue}{gap> XS1 := CrossedSquareByNormalSubgroups( c5d, d10a, d10b, d20 );  }
[  c5d -> d10a ]
[   |  |   ]
[ d10b -> d20  ]
\textcolor{blue}{gap> IsCrossedSquare( XS1 ); }
true
\textcolor{blue}{gap> C2G1 := Cat2GroupOfCrossedSquare( XS1 ); }
(pre-)cat2-group with generating (pre-)cat1-groups:
1 : [((d20 |X d10a) |X (d10b |X c5d))=>(d20 |X d10a)]
2 : [((d20 |X d10a) |X (d10b |X c5d))=>(d20 |X d10b)]
\textcolor{blue}{gap> IsCat2Group( C2G1 ); }
true
\textcolor{blue}{gap> Xab := CrossedSquareOfCat2Group( C2ab ); }
crossed square with crossed modules:
up = [Group( [ (1,5)(2,6)(3,7)(4,8) ] ) -> Group( [ ( 2, 6)( 4, 8) ] )]
left = [Group( [ (1,5)(2,6)(3,7)(4,8) ] ) -> Group(
[ (1,2,3,4)(5,6,7,8), (), () ] )]
right = [Group( [ ( 2, 6)( 4, 8) ] ) -> Group( () )]
down = [Group( [ (1,2,3,4)(5,6,7,8), (), () ] ) -> Group( () )]
\textcolor{blue}{gap> IsCrossedSquare( Xab ); }
true
\textcolor{blue}{gap> IdGroup( Xab ); }
[ [ 2, 1 ], [ 2, 1 ], [ 4, 1 ], [ 1, 1 ] ]
\end{Verbatim}

\section{Table of cat$^{2}$-groups}

The function \textbf{AllCat2Groups(G)} constructs a list $L_G$ 
of all $n_G$ cat$^{2}$-groups $(G;t_1,h_1;t_2,h_2)$ over $G$. 
The function \textbf{AreIsomorphicCat2Groups} is used for checking whether 
or not two cat$^{2}$-groups are isomorphic, 
and \textbf{AllCat2GroupsUpToIsomorphism} returns a list of representatives 
of the isomorphism classes. 
The function \textbf{AllCat2GroupFamilies} returns a list of the positions 
$[1 \ldots n_G]$ partitioned according to these classes. 

In the following \GAP\ session, we compute all $47$ cat$^{2}$-groups on 
$C_{4} \times C_{2}$; representatives of the $14$ isomorphism classes; 
and the list of lists of positions in the families. 
So the eighth class consists of cat$^2$-group numbers $[31,34,35,38]$, 
and \texttt{iso82[8]=all82[31]}. 

\begin{Verbatim}[frame=single, fontsize=\small, commandchars=\\\{\}]
\textcolor{blue}{gap> c4c2 := SmallGroup( 8, 2 );;}
\textcolor{blue}{gap> all82 := AllCat2Groups( c4c2 );;}
\textcolor{blue}{gap> Length( all82 );}
47
\textcolor{blue}{gap> iso82 := AllCat2GroupsUpToIsomorphism( c4c2 );;}
\textcolor{blue}{gap> Length( iso82 );}
14
\textcolor{blue}{gap> AllCat2GroupFamilies( c4c2 );}
[ [ 1 ], [ 2, 5, 8, 11 ], [ 3, 4, 9, 10 ], [ 14, 17, 22, 25 ], 
  [ 15, 16, 23, 24 ], [ 30 ], [ 6, 7, 12, 13 ], [ 31, 34, 35, 38 ], 
  [ 32, 33, 36, 37 ], [ 18, 19, 26, 27 ], [ 20, 21, 28, 29 ], 
  [ 39, 42, 43, 46 ], [ 40, 41, 44, 45 ], [ 47 ] ]
\textcolor{blue}{gap> iso82[8];}
(pre-)cat2-group with generating (pre-)cat1-groups:
1 : [Group( [ f1, f2, f3 ] ) => Group( [ f2, f2 ] )]
2 : [Group( [ f1, f2, f3 ] ) => Group( [ f2, f1 ] )]
\textcolor{blue}{gap> IsomorphismCat2Groups( all82[31], all82[34] ) = fail;}
false
\end{Verbatim}

In the following tables the groups of size at most $30$ are ordered by their
\GAP\ number. 
For each group $G$ we list the number $|IE(G)|$ of idempotent endomorphisms; 
the number $|\calC^1(G)|$ of cat$^1$-groups on $G$, 
followed by the number of their isomorphism classes; 
and then the number $|\calC^2(G)|$ of cat$^2$-groups on $G$, 
and the number of their isomorphism classes. 
The number of isomorphism classes $\calC^1(G)$ of cat$^{1}$-groups 
is given in \cite{alp-wensley-ijac}. 

We may reduce the size of the table by noting the results for cyclic groups. 
When $G=C_{p^{k}}$ is cyclic, with $p$ prime, 
the only idempotent endomorphisms are the identity and zero maps. 
In this case all the cat$^1$-groups have equal tail and head maps, 
and all isomorphism classes are singletons. 
Similarly, when $G = C_{p_1^{k_1}p_2^{k_2} \ldots p_m^{k_m}}$ is cyclic, 
and its order is the product of $m$ distinct primes $p_i$ 
having multiplicities $k_i$, 
there are $2^m$ idempotent endomorphisms and cat$^1$-groups.  
We thus obtain the following results up to $m=4$. 
The rows headed ``groups'' list, for each cat$^2$-group, its four groups 
$[G,R_1,R_2,R_{12}]$ where, for example, $2 \times [G,I,C_{p^k},I]$ 
denotes $\{[G,I,C_{p_1^{k_1}},I],[G,I,C_{p_2^{k_2}},I]\}$. \\ 

\bigskip
\begin{longtable}{ccccccc}
	\hline\hline
	& $G$ 
	    & $|\mathrm{IE}(G)|$ 
	        & $|\calC^1(G)|$ 
	            & $|\calC^1(G)/\cong |$ 
	                & $|\calC^{2}(G)|$ 
	                    & $|\calC^{2}(G)/\cong |$ \\ 
    \hline\hline 
	& $C_{p_1^{k_1}}$ 
	    & 2 
	        & 2 
	            & 2 
	                & 3 
	                    & 3 \\ 
	\hline
	\multicolumn{7}{l}{%
	  \begin{tabular}{ll}
		groups  & $[G,I,I,I],~ [G,I,G,I],~ [G,G,G,G]$  
	  \end{tabular}%
	} \\ 
	\hline\hline 
	& $C_{p_1^{k_1}p_2^{k_2}}$ 
	    & 4 
	        & 4 
	            & 4 
	                & 10 
	                    & 10 \\ 
	\hline
	\multicolumn{7}{l}{ 
	  \begin{tabular}{ll}
		groups & $[G,I,I,I],~ [G,I,G,I],~ [G,G,G,G],~ 2 \times [G,I,C_{p^k},I],$ \\
			   & $2 \times [G,C_{p^k},G,C_{p^k}],~ 
			      2 \times [G,C_{p^k},C_{p^k},C_{p^k}],~ 
				 [G,C_{p_1^{k_1}},C_{p_2^{k_2}},I]$ 
	  \end{tabular} 
	} \\ 
	\hline\hline 
	& $C_{p_1^{k_1}p_2^{k_2}p_3^{k_3}}$ 
	    & 8 
	        & 8 
	            & 8 
	                & 36 
	                    & 36 \\ 
	\hline
	\multicolumn{7}{l}{ 
	  \begin{tabular}{ll}
		groups & $[G,I,I,I],~ [G,I,G,I],~ [G,G,G,G],~ 3 \times [G,I,C_{p^k},I],$ \\
			   & $3 \times [G,C_{p^k},G,C_{p^k}],~ 
			      3 \times [G,C_{p^k},C_{p^k},C_{p^k}],~ 
				  3 \times [G,C_{p^k},C_{q^j},I],$ \\
			   & $3 \times [G,I,C_{p^kq^j},I],~ 
			      3 \times [G,C_{p^kq^j},G,C_{p^kq^j}],~ 
				  3 \times [G,C_{p^kq^j},C_{p^kq^j},C_{p^kq^j}],$ \\ 
			   & $6 \times [G,C_{p^k},C_{p^kq^j},C_{p^k}],~ 
			      3 \times [G,C_{r^i},C_{p^kq^j},I],~ 
				  3 \times [G,C_{p^kq^j},C_{p^kr^i},C_{p^k}]$ 
		\end{tabular} 
	} \\ 
	\hline\hline 
	& $C_{p_1^{k_1}p_2^{k_2}p_3^{k_3}p_4^{k_4}}$ 
	    & 16 
	        & 16 
	            & 16 
	                & 136 
	                    & 136 \\ 
	\hline\hline
\end{longtable}

\noindent 
When $m=1$ there are $16$ cyclic groups of order at most $30$; 
when $m=2$ there are $12$ such groups; 
and when $m=3$ there is just the small group $30/4 = C_{30}$.

The following table contains the results for those $G$ which are not cyclic. 
\begin{longtable}{ccrrrrr}
	\hline 
	{\GAP\ }\# 
	    & $G$ 
	        & $|\mathrm{IE}(G)|$ 
	            & $|\calC^1(G)|$ 
	                & $|\calC^1(G)/\cong |$ 
	                    & $|\calC^{2}(G)|$ 
	                        & $|\calC^{2}(G)/\cong |$  \\ 
	\hline
	1/1 & $I$ & 1 & 1 & 1 & 1 & 1 \\ 
	4/2 & $K_4 = C_2 \times C_2$ & 8 & 14 & 4 & 36 & 9 \\ 
	6/1 & $S_3$ & 5 & 4 & 2 & 7 & 3 \\ 
	8/2 & $C_4 \times C_2$ & 10 & 18 & 6 & 47 & 14 \\ 
	8/3 & $D_8$ & 10 & 9 & 3 & 21 & 6 \\ 
	8/4 & $Q_8$ & 2 & 1 & 1 & 1 & 1 \\ 
	8/5 & $C_2 \times C_2 \times C_2$ & 58 & 226 & 6 & 1,711 & 23 \\ 
	9/2 & $C_3 \times C_3$ & 14 & 38 & 4 & 93 & 9 \\ 
	10/1 & $D_{10}$ & 7 & 6 & 2 & 11 & 3 \\ 
	12/1 & $C_3 \ltimes C_4$ & 5 & 4 & 2 & 7 & 3 \\ 
	12/3 & $A_4$ & 6 & 5 & 2 & 9 & 3 \\ 
	12/4 & $D_{12}$ & 21 & 12 & 4 & 41 & 10 \\ 
	12/5 & $C_3 \times K_4$ & 16 & 28 & 8 & 136 & 32 \\ 
	14/1 & $D_{14}$ & 9 & 8 & 2 & 15 & 3 \\ 
	16/2 & $C_4 \times C_4$ & 26 & 98 & 5 & 231 & 11 \\ 
	16/3 & $(C_4 \times C_2) \ltimes C_2$ & 18 & 25 & 4 & 57 & 7 \\ 
	16/4 & $C_4 \ltimes C_4$ & 10 & 17 & 3 & 25 & 4 \\ 
	16/5 & $C_8 \times C_2$ & 10 & 18 & 6 & 47 & 14 \\ 
	16/6 & $C_8 \ltimes C_2$ & 6 & 5 & 2 & 9 & 3 \\ 
	16/7 & $D_{16}$ & 18 & 9 & 2 & 17 & 3 \\ 
	16/8 & $QD_{16}$ & 10 & 5 & 2 & 9 & 3 \\ 
	16/9 & $Q_{16}$ & 2 & 1 & 1 & 1 & 1 \\ 
	16/10 & $C_4 \times K_4$ & 82 & 322 & 12 & 2,875 & 53 \\ 
	16/11 & $C_2 \times D_8$ & 82 & 97 & 9 & 649 & 29 \\ 
	16/12 & $C_2 \times Q_8$ & 18 & 17 & 3 & 25 & 4 \\ 
	16/13 & $(C4 \times C2) \ltimes C_2$ & 26 & 13 & 2 & 37 & 4 \\ 
	16/14 & $K_4 \times K_4$ & 382 & 4,162 & 9 & 298,483 & 53  \\ 
	18/1 & $D_{18}$ & 11 & 10 & 2 & 19 & 3 \\ 
	18/3 & $C_3 \times S_3$ & 12 & 8 & 4 & 24 & 10 \\ 
	18/4 & $(C_3 \times C_3) \ltimes C_2$ & 47 & 118 & 4 & 541 & 9 \\ 
	18/5 & $C_6 \times C_3$ & 28 & 76 & 8 & 358 & 32 \\ 
	20/1 & $Q_{20}$ & 7 & 6 & 2 & 11 & 3 \\ 
	20/3 & $C_4 \ltimes C_5$ & 7 & 6 & 2 & 11 & 3 \\ 
	20/4 & $D_{20}$ & 31 & 18 & 4 & 65 & 10 \\ 
	20/5 & $C_5 \times K_4$ & 16 & 28 & 8 & 136 & 32 \\ 
	21/1 & $C_3 \ltimes C_7$ & 9 & 8 & 2 & 15 & 3 \\ 
	22/1 & $D_{22}$ & 13 & 12 & 2 & 23 & 3 \\ 
	24/1 & $C_3 \ltimes C_8$ & 5 & 4 & 2 & 7 & 3 \\ 
	24/3 & $SL(2,3)$ & 6 & 1 & 1 & 1 & 1 \\ 
	24/4 & $Q_{24}$ & 5 & 4 & 2 & 7 & 3 \\ 
	24/5 & $S_3 \times C_4$ & 27 & 12 & 4 & 41 & 10 \\ 
	24/6 & $D_{24}$ & 33 & 20 & 4 & 75 & 10 \\ 
	24/7 & $Q_{12} \times C_2$ & 25 & 36 & 6 & 115 & 14 \\ 
	24/8 & $D_8 \ltimes C_3$ & 23 & 12 & 4 & 41 & 10 \\ 
	24/9 & $C_{12} \times C_2$ & 20 & 36 & 12 & 178 & 52 \\ 
	24/10 & $D_8 \times C_3$ & 20 & 18 & 6 & 75 & 20 \\ 
	24/11 & $Q_8 \times C_3$ & 4 & 2 & 2 & 3 & 3 \\ 
	24/12 & $S_4$ & 12 & 5 & 2 & 9 & 3 \\ 
	24/13 & $A_4 \times C_2$ & 15 & 10 & 4 & 31 & 10 \\ 
	24/14 & $S_3 \times K_4$ & 157 & 116 & 8 & 999 & 32 \\ 
	24/15 & $C_6 \ltimes K_4$ & 116 & 452 & 12 & 6,786 & 84 \\ 
	25/2 & $C_5 \times C_5$ & 32 & 152 & 4 & 348 & 9 \\ 
	26/1 & $D_{26}$ & 15 & 14 & 2 & 27 & 3 \\ 
	27/2 & $C_9 \times C_3$ & 20 & 56 & 6 & 138 & 14 \\ 
	27/3 & $(C_3 \times C_3) \ltimes C_3$ & 38 & 37 & 2 & 127 & 4 \\ 
	27/4 & $C_9 \ltimes C_3$ & 11 & 10 & 2 & 19 & 3 \\ 
	27/5 & $C_3 \times C_3 \times C_3$ & 236 & 2,108 & 6 & 24,222  & 16 \\ 
	28/1 & $Q_{28}$ & 9 & 8 & 2 & 15 & 3 \\ 
	28/3 & $D_{28}$ & 41 & 24 & 4 & 89 & 10 \\ 
	28/4 & $C_7 \times K_4$ & 16 & 28 & 8 & 136 & 32 \\ 
	30/1 & $S_3 \times C_5$ & 10 & 8 & 4 & 24 & 10 \\ 
	30/2 & $D_{10} \times C_3$ & 14 & 12 & 4 & 38 & 10 \\ 
	30/3 & $D_{30}$ & 25 & 24 & 4 & 92 & 10 \\ 
	\hline
\end{longtable}
The $1,000$ isomorphism classes contain just $13$ cat$^2$-groups 
whose diagonal is \emph{not} a cat$^1$-group: 
one each for groups [8/3, 16/3, 16/13, 27/3], 
three for 24/10 and six for 16/11. 

\section*{Acknowledgement}

The first and second authors were supported by Eskisehir Osmangazi
University Scientific Research Projects (Grant No: 2017/19033).

\end{document}